\documentclass[12pt,sort&compress]{article}
\usepackage{amsfonts}
\usepackage{bm}
\usepackage{makecell}
\usepackage{amssymb,amsmath}
\usepackage{amscd}
\usepackage{lipsum}
\usepackage{latexsym}
\usepackage{CJK}
\usepackage{longtable}
\usepackage{color}
\usepackage{enumerate}
\def\a{\alpha}
\def\g{\gamma}

\def\b{\beta}
\def\cl{\centerline}

\def\vs{\vspace*}

\def\Z{\mathbb{Z}}
\def\C{\mathbb{C}}
\def\QED{\hfill$\Box$}

\def\pa{\partial}
\def\la{\lambda}

\def\HV{\mathcal{HV}}

\def\A{\mathcal{A}}

\numberwithin{equation}{section}
\newtheorem{theo}{Theorem}[section]
\newtheorem{defi}[theo]{Definition}
\newtheorem{coro}[theo]{Corollary}
\newtheorem{lemm}[theo]{Lemma}
\newtheorem{prop}[theo]{Proposition}

\newtheorem{rema}[theo]{Remark}

\newtheorem{example}[theo]{Example}
\newtheorem{proposition}[theo]{Proposition}
\makeatletter
\newcommand\blfootnote[1]{%
\begingroup
\renewcommand\thefootnote{}\footnote{#1}%
\addtocounter{footnote}{-1}%
\endgroup
}

\def\@biblabel#1{#1.~}

\makeatother

\begin{document}
\vs{10pt} \cl{\large {\bf A class of graded conformal algebras which is induced by }}
\cl{\large {\bf Heisenberg-Virasoro conformal algebra}}
\vs{10pt}\cl{Lipeng Luo$^1$, Yucai Su$^2$ and Xiaoqing Yue$^3$}
\cl{\small{School of Mathematical Sciences, Tongji University, Shanghai, 200092, China.}}
\cl{\small E-mail: luolipeng@tongji.edu.cn, ycsu@tongji.edu.cn, xiaoqingyue@tongji.edu.cn.
 }\vs{6pt}
\small\parskip .005 truein \baselineskip 3pt \lineskip 3pt

\noindent{\bf Abstract:} In this paper, we obtain a class of $\Z$-graded conformal algebras which is induced by Heisenberg-Virasoro conformal algebra. More precisely, we classify $\Z$-graded conformal algebras ${\A}=\oplus^\infty_{i=-1}\A_i$ satisfying the following conditions,\vs{5pt}\\
(C1) $\A_0$ is the Heisenberg-Virasoro conformal algebra;\\
(C2) Each $\A_i$ for $i\in\Z_{\ge-1}^*$ is an $\A_0$-module of rank one;\\
(C3) $[{X_{-1}}_\lambda X_i]\neq 0$ for $i\ge 0$, where $X_i$ is any one of $\C[\partial]$-generators of $\A_i$ for $i\in \Z_{\ge -1}$.\vs{5pt}\\

\noindent Further, we prove that all finite nontrivial irreducible modules of these algebras under some special conditions are free of rank one as a $\C[\partial]$-module. The conformal derivations of this class of graded Lie conformal algebras are also determined.
\blfootnote{Corresponding author: Xiaoqing Yue (xiaoqingyue@tongji.edu.cn)}

\vs{5pt}

\vs{5pt}


\section{Introduction}
Lie conformal algebra, first introduced by Kac in \cite{K1,K2}, gives an axiomatic description of the operator product expansion of chiral fields in conformal field theory (see \cite{BPZ}). It has been shown that the theory of Lie conformal algebras has close connections to vertex algebras, infinite-dimensional Lie algebras satisfying the locality property in \cite{KacL} such as affine Kac-Moody algebras, the Virasoro algebra, and the Hamiltonian formalism in the theory of nonlinear evolution equations (see \cite{BDK,DK,K2,SYX, SY,Z1,Z2}).

In recent years, the general structure theory, the representation theory and the cohomology theory of Lie conformal algebras were systematically developed in \cite{CK,DK,BKV}. The Virasoro Lie conformal algebra $Vir$ and the current Lie conformal algebra $Cur\mathcal{G}$ associated to a Lie algebra $\mathcal{G}$  are two important examples of Lie conformal algebras. It was shown in \cite{DK} that $Vir$ and all current Lie conformal algebras $Cur\mathcal{G}$ where $\mathcal{G}$ is a finite dimensional simple Lie algebra exhaust all finite simple Lie conformal algebras. Finite irreducible conformal modules of these simple Lie conformal algebras were classified in \cite{CK} by investigating the representation theory of their extended annihilation algebras. The general cohomology theory of conformal algebras with coefficients in an arbitrary conformal module was
developed in \cite{BKV}, where explicit computations of cohomology groups for the Virasoro conformal algebra and current conformal algebras were given. Some low dimensional cohomology groups of the general Lie
conformal algebras $gc_N$ were studied in \cite{S}.
All cohomology groups of the Heisenberg-Virasoro conformal algebra with trivial coefficients were determined in \cite{YW}.
However, the general structure theory, the representation theory and the cohomology theory of infinite Lie conformal algebras are very difficult to investigate.

For infinite Lie conformal algebras, the general Lie conformal algebra $gc_N$ plays an important role in the theory of Lie conformal algebras as the general Lie algebra $gl_N$ does in the theory of Lie algebras. So there are lots of papers about the general Lie conformal algebra $gc_N$ revealed in \cite{BKL1,BKL2,GXY,S,SXY,SYX,WCY}. In particular, in \cite{SYX2}, Su and Yue obtain a simple $\Z$-graded Lie conformal algebras $\mathcal{G}=\oplus^\infty_{i=-1}\mathcal{G}_i$, which satisfies that $\mathcal{G}_0$ is the Virasoro conformal algebra and each $\mathcal{G}_i$ for $i\in\Z_{\ge-1}$ is a $Vir$-module of rank one. Moreover, they declare that such a Lie conformal algebra $\mathcal{G}$, which is a finitely freely generated simple Lie conformal algebra of linear growth that cannot be embedded into $gc_N$ for any $N$, does not have a nontrivial representation on any finite $\C[\pa]$-module. In \cite{XH}, Xu and Hong classify all finite irreducible modules of a class of $\Z_+$-graded Lie conformal algebras $\mathcal{L}=\oplus^\infty_{i=0}\C[\pa]L_i$ satisfying $[{L_0}_\la {L_0}]=(\pa+2\la)L_0$ and $[{L_1}_\la L_i]\ne0$ for any $i\in\Z_+$. Furthermore, they show that all finite nontrivial irreducible modules of these algebras are free of rank one as a $\C[\pa]$-module. Inspired by the above results, in this paper, we aim to investigate the structure theory and the representation theory of all graded Lie conformal algebras satisfying (C1)--(C3). This is the main motivation to present our work.

This paper is organized as follows. In Section 2, we recall some basic definitions, notations, and related known results about Lie conformal algebras. In Section 3, we classify all $\Z$-graded conformal algebras ${\A}$ and study the conformal algebra structure of ${\A}$, which satisfies conditions (C1)--(C3). In Section 4, we investigate some properties of graded Lie conformal algebras ${\HV}(\a,\b)$ defined in Definition \ref{def2.7}. Moreover, we classify all finite nontrivial irreducible modules over graded Lie conformal algebras ${\HV}(\a,\b)$ with $\a\ne1$ and $\b\ne0$. In Section 5, we show that all conformal derivations of ${\HV}(\a,\b)$ are inner.
Our main results are summarized in Theorems \ref{th3.5}, \ref{th4.11} and \ref{th5.7}.

Throughout this paper, we use notations $\mathbb{C}$, $\mathbb{C}^*$, $\mathbb{Z}$, $\mathbb{Z_{+}}$, $\mathbb{Z}_{\ge-1}$, $\mathbb{Z}_{\ge-1}^*$ and $\mathbb{Z}_{\ge2}$ to represent the set of complex numbers, nonzero complex numbers, integers, nonnegative integers, integers greater than $-2$, nonzero integers greater than $-2$ and integers greater than $1$, respectively. In addition, all vector spaces and tensor products are over $\mathbb{C}$. For convenience, we abbreviate $\otimes_{\mathbb{C}}$ to $\otimes$.
\vs{8pt}

\section{Preliminaries}

   In this section, we first recall some basic definitions, notations and related results about Lie conformal algebras for later use. For a detailed description, one can refer to \cite{BKV,YW}.

\begin{defi}
\begin{em}
A \emph{Lie conformal algebra} $\mathcal {A}$ is a $\C[\partial ]$-module endowed with a $\C$-linear map $\mathcal {A}\otimes \mathcal {A}\rightarrow \C[\lambda]\otimes \mathcal {A}$, $a\otimes b \mapsto [a_\lambda b],
$ and
satisfying the following axioms ($a, b, c\in \mathcal {A}$),
\begin{align}
[\partial a_\lambda b]&=-\lambda[a_\lambda b],\ \ [ a_\lambda \partial b]=(\partial+\lambda)[a_\lambda b] \ \ \mbox{(conformal\  sesquilinearity)},\label{Lc1}\\
{[a_\lambda b]} &= -[b_{-\lambda-\partial}a] \ \ \mbox{(skew-symmetry)},\label{Lc2}\\
{[a_\lambda[b_\mu c]]}&=[[a_\lambda b]_{\lambda+\mu
}c]+[b_\mu[a_\lambda c]]\ \ \mbox{(Jacobi \ identity)}\label{Lc3}.
\end{align}
\end{em}
\end{defi}

The notions of subalgebras, ideals, quotients and homomorphisms of Lie conformal algebras are obvious. Moreover, due to the skew\text{-}symmetry, any left or right ideal is actually a two-side ideal. A Lie conformal algebra $\mathcal{A}$ is called \emph {finite} if $\mathcal{A}$ is finitely generated as a $\mathbb{C}[\partial]$-module. The \emph {rank} of a Lie conformal algebra $\mathcal{A}$, denoted by rank($\mathcal{A}$), is its rank as a $\mathbb{C}[\partial]$-module. A Lie conformal algebra $\mathcal{A}$ is \emph{simple} if it has no nontrivial ideals and it is not abelian.

The Virasoro conformal algebra $Vir$ and the current Lie conformal algebra $Cur\mathcal{G}$ associated to a Lie algebra $\mathcal{G}$  are two important examples of Lie conformal algebras.

\begin{example}
	\em{
		The conformal algebra $Vir$ is a free $\mathbb{C}[\partial]$-module on the generator $L$, which is defined as follows
		\begin{align*}
		Vir= \mathbb{C}[\partial]L, \quad [L_\lambda L]= (\partial+ 2\lambda)L.
		\end{align*}}
\end{example}
\begin{example}
	\em{
		Let $\mathcal{G}$ be a Lie algebra. The current Lie conformal algebra associated to $\mathcal{G}$ is defined as follows
		\begin{align*}
		Cur\mathcal{G}= \mathbb{C}[\partial] \otimes \mathcal{G}, \quad [a_\lambda b]= [a, b],
		\end{align*}
		for all $a, b \in \mathcal{G}$.}
\end{example}
\begin{example}
	\em{
		Let $\mathcal{G}$ be a Lie algebra. The \emph{standard semi-direct sum} of $Vir$ and $Cur\mathcal{G}$ denoted by $Vir\ltimes Cur\mathcal{G}=Vir\oplus Cur\mathcal{G}$ can be given a conformal algebra structure by
		\begin{align*}
		[L_\la L]= (\pa+ 2\la)L,\quad[L_\la a]= (\pa+ \la)a,\quad [a_\la b]= [a, b],
		\end{align*}
		for all $a, b \in \mathcal{G}$.}
\end{example}

It was shown in \cite{DK} that $Vir$ and all current Lie conformal algebras $Cur\mathcal{G}$ where $\mathcal{G}$ is a finite dimensional simple Lie algebra exhaust all finite simple Lie conformal algebras.

Let $\A$ be a Lie conformal algebra. There is an important Lie algebra associated to it. For each $j \in\mathbb{Z_{+}}$, we can define the  \emph {$j$-th product} $a_{(j)}b$ of two elements $a,b \in\A$ as follows:
\begin{align}
  [a_\lambda b]=\sum_{j\in\mathbb{Z_{+}}}(a_{(j)}b)\frac{\lambda^{j}}{j!}.
\end{align}
Let $Lie(\A)$ be the quotient
of the vector space with basis $a_{(n)}$ $(a\in \A, n\in\mathbb{Z})$ by
the subspace spanned over $\mathbb{C}$ by
elements:
$$(\alpha a)_{(n)}-\alpha a_{(n)},~~(a+b)_{(n)}-a_{(n)}-b_{(n)},~~(\partial
a)_{(n)}+na_{(n-1)},$$ where $a,b\in \A, \alpha\in \mathbb{C}, n\in
\mathbb{Z}$. The operation on $Lie(\A)$ is defined as follows:
\begin{equation}\label{106}
[a_{(m)}, b_{(n)}]=\sum_{j\in \mathbb{Z_{+}}}\left(\begin{array}{ccc}
m\\j\end{array}\right)(a_{(j)}b)_{(m+n-j)}.\end{equation} Then
$Lie(\A)$ is a Lie algebra and it is called the\emph{ coefficient algebra} of $\A$ (see \cite{K2}).

\begin{defi}
\begin{em}
The \emph {annihilation algebra} associated to a Lie conformal algebra $\A$ is the subalgebra
\begin{align}
  Lie(\A)^{+}=span_\C\{a_{(n)}\ |\ a\in \A, n \in \mathbb{Z_{+}}\},\nonumber
\end{align}
 of the Lie algebra $Lie(\A)$. The semi-direct sum of the 1-dimensional Lie algebra $\mathbb{C}\partial$ and $Lie(\A)^{+}$ with the action $ \partial (a_{(n)})=-na_{(n-1)}$ is called the \emph {extended annihilation algebra} $Lie(\A)^e$.
\end{em}		
\end{defi}

The derived algebra of a Lie conformal algebra $\A$ is the vector space $\A^{'}=span_{\C}\{a_{(i)}b\ |\ a,b\in\A,\ i\in\Z_+\}$. It is not difficult to check that $\A^{'}$ is an ideal of $\A$. Define $\A^{(1)}=\A^{'}$ and $\A^{(n+1)}=\A^{(n)'}$ for any $n\ge1$. A Lie conformal algebra $\A$ is said to be \emph{solvable}, if there exists some $n\in\Z_+$ such that $\A^{(n)}=0$. Obviously, the sum of two solvable ideals of a Lie conformal algebra is still a solvable ideal. Hence, if
$\A$ is a finite Lie conformal algebra, the maximal solvable ideal exists and is unique, which is called \emph{radical} of $\A$ and denoted by $Rad(\A)$.
\begin{defi}
\begin{em}
A Lie conformal algebra is \emph{semisimple} if its radical is zero. In particular, $\A/Rad(\A)$ is semisimple for any finite Lie conformal algebra $\A$.
\end{em}
\end{defi}

It is revealed in Theorem 7.1 in \cite{DK} that any finite semisimple Lie conformal algebra is the direct sum of the following Lie conformal algebras:
$$Vir,\  Cur\mathcal{G},\  Vir\ltimes Cur\mathcal{G},$$ where $\mathcal{G}$ is a finite dimensional semisimple Lie algebra.

There is an important example of Lie conformal algebra of rank two, called the Heisenberg-Virasoro conformal algebra and denoted by $\mathcal{HV}$, which is not a simple Lie  conformal algebra, as shown in \cite{YW}.

\begin{example}
	\em{
		The Lie conformal algebra $\mathcal{HV}= \mathbb{C}[\partial]L\oplus \mathbb{C}[\partial]H$ is a free $\mathbb{C}[\partial]$-module on the generator $L, H$, which is defined as follows
\begin{align}\label{CHV}
[L_\lambda L]= (\partial+ 2\lambda)L, \quad [L_\lambda H]= (\partial+ \lambda)H, \quad [H_\lambda L]= \lambda H, \quad [H_\lambda H]=0.
\end{align}}
\end{example}
It is not difficult to see that $\mathcal{HV}$ contains the Virasoro conformal algebra $Vir$ as a subalgebra. Moreover, $\mathcal{HV}$ has a nontrivial abelian conformal ideal with one free generator $H$ as a $\mathbb{C}[\partial]$-module. Obviously, it is neither simple nor semisimple.

 Let us introduce a class of graded conformal algebras ${\HV}(\a,\b)$, which is closely related to the Heisenberg-Virasoro conformal algebra.
 \begin{defi}\label{def2.7}
 \begin{em}The
Lie conformal algebra ${\HV}(\a,\b)=\oplus^\infty_{i=-1}{\HV}(\a,\b)_i$ is a graded conformal algebra with $\C[\pa]$-basis $\{L, H_i\ |\ i\in\Z_{\ge -1}\}$, where $H_i$ is a $\C[\partial]$-generator of ${\HV}(\a,\b)_i$ for $i\in \Z_{\ge -1}$ and ${\HV}(\a,\b)_0=\C[\pa]L\oplus\C[\pa]H_0\cong\HV$ such that
 \begin{align}\label{eq2.5}
[L_\la L]=(\pa+2\la)L, \quad [L_\la H_i]=\big(\pa+(i\a-i+1)\la+i\b\big)H_i, \quad [{H_i}_\la H_j]=(j-i)H_{i+j},
 \end{align}
for $i,j\in\Z_{\ge-1}$ and $\a,\b\in\C$.
 \end{em}
 \end{defi}

The \emph{general Lie conformal algebra $gc_N$} can be defined as the infinite rank $\C[\pa]$-module $\C[\pa,x]\otimes gl_N$, with the $\lambda$-bracket
\begin{align}
[f(\pa,x)A_{\,\lambda\,} g(\pa,x)B]=f(-\lambda,x+\pa+\la)g(\pa+\la,x)AB-f(-\la,x)g(\pa+\la,x-\la)BA,
\end{align}
for $f(\pa,x),g(\pa,x)\in\C[\pa,x],$\ $A,B\in gl_N$, where $gl_N$ is the space of $N\times N$ matrices, and we have identified $f(\pa,x)\otimes A$ with $f(\pa,x)A$. If we set $J_A^n=x^nA$, then
\begin{align}
[{J_A^m}_\la J_B^n]=\sum^m_{s=0}\binom m s (\la+\pa)^s J_{AB}^{m+n-s}-\sum^n_{s=0}\binom n s (-\la)^s J_{BA}^{m+n-s},
\end{align}
for $m,n\in\Z_+,$\ $A,B\in gl_N,$ where $\binom m s =\frac{m(m-1)\cdots(m-s+1)}{s!}$ if $s\ge0$ and $\binom m s=0$ otherwise, is the binomial coefficient.

\begin{defi}
\begin{em}
Let $\A$ be a Lie conformal algebra. For any $x\in\A$, we define the operator $(ad\ x)_\la:\A\to\A[\la]$ such that $(ad\ x)_\la(y)=[x_\la y]$ for any $y\in\A$. An element $x\in\A$ is \emph{locally nilpotent} if for any $y\in\A$, there exists $n\in\Z_{\ge1}$ such that $(ad\ x)_{\la}^n(y)=0$.

\end{em}
\end{defi}

\begin{defi}
\begin{em}
A \emph{conformal module} $V$ over a Lie conformal algebra $\mathcal {A}$
is a $\mathbb{C}[\partial]$-module endowed with a $\C$-linear map
$\mathcal {A}\otimes V\rightarrow
V[\lambda]$, $a\otimes v\mapsto a_\lambda v$, such that for any $a,b\in\mathcal {A}$, $v\in V$,
\begin{eqnarray*}
&&a_\lambda(b_\mu v)-b_\mu(a_\lambda v)=[a_\lambda b]_{\lambda+\mu}v,\\
&&(\partial a)_\lambda v=-\lambda a_\lambda v,\ a_\lambda(\partial
v)=(\partial+\lambda)a_\lambda v.
\end{eqnarray*}
 If $V$ is finitely generated over $\mathbb{C}[\partial]$, then $V$ is simply called \emph {finite}. The \emph {rank} of a conformal module $V$ is its rank as a $\mathbb{C}[\partial]$-module. A conformal module $V$ is called \emph {irreducible} if it has no nontrivial submodules.
\end{em}
\end{defi}

In the following, since we only consider conformal modules, we abbreviate ``conformal modules" to ``modules".

Let $\A$ be a Lie conformal algebra and $V$ be an $\A$-module. An element $v\in V$ is called \emph {invariant} if $\A_\lambda v=0$. Obviously, the set of all invariants of $V$ is a conformal submodule of $V$, denoted by $V^0$. An $\A$-module $V$ is called trivial if $V^0=V$, i.e., a module on which $\A$ acts trivially. The vector space $\C$ can be regarded as a trivial module with trivial actions of $\partial$ and $\mathcal {A}$. For a fixed nonzero complex constant $a$, there
is a natural $\C[\partial]$-module $\C_a$, such that $\C_a=\C$ and $\partial v=a v$ for $v\in\C_a$. Then
$\C_a$ becomes an $\mathcal{A}$-module with
$\mathcal{A}$ acting by zero. It is easy to check that the modules $\mathbb{C}_a$ with $a \in \mathbb{C}$ exhaust all trivial irreducible $\mathcal{A}$-modules.

\begin{defi}\label{def2.11}
\begin{em}
Let $\A$ be a Lie conformal algebra and $V$ be an $\A$-module. If $a_\la V=0$ for some $a\in\A$ if and only if $a=0$, then $V$ is called a \emph{faithful} $\A$-module.
\end{em}
\end{defi}

For any $\A$-module $V$, we have some basic results as follows.

\begin{lemm}\label{lem1}\begin{em}(\cite{KacL}, Lemma 2.2)\end{em}
Let $\A$ be a Lie conformal algebra and $V$ be an $\A$-module.
\begin{enumerate}
 \item [\rm(1)] If $\partial v=av$ for some $a \in \mathbb{C}$ and $v \in V$, then $\A_\lambda v=0$.
 \item [\rm(2)] If $V$ is a finite module without any nonzero invariant element, then $V$ is a free $\mathbb{C}[\partial]$-module.
\end{enumerate}
\end{lemm}

 Let $V$ be an $\A$-module. An element $v \in V$ is called a \emph {torsion element} if there exists a nonzero polynomial $ p(\partial)\in \mathbb{C}[\partial]$ such that $p(\partial)v=0$. For any $\mathbb{C}[\partial]$-module $V$, it is not difficult to check that there exists a nonzero torsion element if and only if there exists nonzero $v \in V$ such that $\partial v=av$ for some $a \in \mathbb{C}$ by using \emph {The Fundamental Theorem of Algebra}. By Lemma \ref{lem1}, we can deduce that a finitely generated $\mathbb{C}[\partial]$-module is free if and only if it has no nonzero torsion element. Thus, the following result is obvious.

 \begin{lemm}
Let $\A$ be a Lie conformal algebra and $V$ be a finite nontrivial irreducible $\A$-module. Then $V$ has no nonzero torsion elements and is free of a finite rank as a $\mathbb{C}[\partial]$-module.
\end{lemm}

Similar to the definition of the $j$-th product $a_{(j)}b$ of two elements $a,b \in\A$, we can also define \emph {$j$-th actions} of $\A$ on $V$ for each $j \in\mathbb{Z_{+}}$, i.e., $a_{(j)}v$ for any $a \in\A,v\in V$ by
\begin{align}
  a_\lambda v=\sum_{j\in\mathbb{Z_{+}}}(a_{(j)}v)\frac{\lambda^{j}}{j!}.
\end{align}
In \cite{CK}, Cheng and Kac investigated a close connection between the module of a Lie conformal algebra and that of its extended annihilation algebra.

\begin{lemm}\label{ll1}
Let $\A$ be a Lie conformal algebra and $V$ be an $\A$-module. Then $V$ is precisely a module over $Lie(\A)^e$ satisfying the property
\begin{align}\label{eq1}
  a_{(n)}\cdot v=0,\quad n\geq N,
\end{align}
for $a \in \A, v \in V$, where $N$ is a non-negative integer depending on $a$ and $v$.
\end{lemm}

\begin{rema}
\begin{em}
By abuse of notations, we also call a Lie algebra module satisfying (\ref{eq1}) a conformal module over $Lie(\mathcal{A})^e$.
\end{em}
\end{rema}

It is well known from \cite{CK, WY} that

 \begin{proposition}\label{pr1}
 All free nontrivial $Vir$-modules of rank one over $\mathbb{C}[\partial]$ are as follows \begin{em}($ \alpha,\beta \in \mathbb{C}$):\end{em}
  \begin{align*}
V_{\alpha,\beta}=\mathbb{C}[\partial]v,\qquad L_\lambda v=(\partial+\alpha\lambda+\beta)v.
 \end{align*}
 Moreover, the module $V_{\alpha,\beta}$ is irreducible if and only if $\alpha$ is non-zero. The module $V_{0,\beta}$ contains a unique nontrivial submodule $ (\partial+\beta)V_{0,\beta}$ isomorphic to $V_{1,\beta}$. The modules $V_{\alpha,\beta}$ with $\alpha\neq 0$ exhaust all finite irreducible nontrivial $Vir$-modules.
 \end{proposition}

  \begin{proposition}\label{pr2}
 All free nontrivial $\mathcal{HV}$-modules of rank one over $\mathbb{C}[\partial]$ are as follows \begin{em}($ \alpha,\beta,\gamma \in \mathbb{C}$):\end{em}
  \begin{align}\label{MHV}
V_{\alpha,\beta, \gamma}=\mathbb{C}[\partial]v,\qquad L_\lambda v=(\partial+\alpha\lambda+\beta)v, \quad H_\lambda v=\gamma v.
 \end{align}
 In particular, $V_{\alpha,\beta, \gamma}$ is an irreducible ${\HV}$-module if and only if $\a\ne0$ or $\gamma \ne0$.
 \end{proposition}

 \begin{rema}
 It is revealed in \begin{em}\cite{CK, WY}\end{em} that all finite nontrivial irreducible $Vir$(or $\mathcal{HV}$)-modules are free of rank one.
 \end{rema}


\section{Classification of graded Lie conformal algebras}
In this section, we will classify all $\Z$-graded Lie conformal algebras ${\A}$ and study the Lie conformal algebra structure of ${\A}$, which satisfies conditions (C1)--(C3).

Let ${\A}=\oplus^\infty_{i=-1}\A_i$ be a Lie conformal algebra satisfying conditions (C1)--(C3) and denote by $X_i$ a $\C[\partial]$-generator of $\A_i$ for $i\in \Z_{\ge -1}^*$. For convenience, we always use $L$ and $X_0(:=H)$ to denote the $\C[\partial]$-generators of $\A_0(:=\HV)$ as shown in (\ref{CHV}) in the sequel.

By the assumption of ${\A}$ and (\ref{MHV}), we can suppose that
\begin{align}
&[L_\la X_i]=(\pa+\a_i\la+\b_i)X_i,\quad [H_\la X_i]=\g_iX_i,\quad for\  i\in\Z_{\ge-1},\label{eq3.1}\\
&[{X_i}_\la X_j]=f_{i,j}(\pa,\la)X_{i+j},\quad for\  i,j\neq 0\  and \ i+j\neq 0,\label{eq3.2}\\
&[{X_{-1}}_\la X_1]=g_{-1,1}(\pa,\la)L+f_{-1,1}(\pa,\la)H,\quad [{X_{-1}}_\la X_{-1}]=0,\label{eq3.3}
\end{align}
where $\a_i,\b_i,\g_i\in\C$ and $g_{-1,1}(\pa,\la),f_{i,j}(\pa,\la)\in \C[\pa,\la]$ are polynomials of $\pa$ and $\la$ for $i,j\in\Z_{\ge-1}$. In particular, we can regard $\g_j$ as $f_{0,j}(\pa,\la)$ in (\ref{eq3.2}) for $j\in\Z_{\ge-1}$. Because of (C3), it is not difficult to see that $X_{i+j}$ can be generated by $X_i$ and $X_j$ for $i,j\in\Z_{\ge-1}$.

To simplify formula (\ref{eq3.1}), we can introduce the following lemma.
\begin{lemm}\label{lm3.1}
$\b_i=i\b_1$ and $\g_i=i\g_1$ for $i\in\Z_{\ge-1}$, $\b_1\in\C$, $\g_1\in\C^*$. Moreover,
\begin{align}
&[L_\la X_i]=(\pa+\a_i\la+i\b_1)X_i,\quad [H_\la X_i]=i\g_1X_i,\quad for\  i\in\Z_{\ge-1}.
\end{align}
\end{lemm}
{\it Proof.} Using the Jacobi identity $[L_{\mu}[{X_{-1}}_\la X_i]]=[[L_{\mu}X_{-1}]_{\la+\mu} X_i]+[{X_{-1}}_\la [L_{\mu} X_i]]$ and comparing the coefficients of $X_{i-1}$ for $i\in\Z_{\ge2}$, we can obtain that
\begin{align}\label{eq3.5}
(\pa+&\a_{i-1}\mu+\b_{i-1})f_{-1,i}(\pa+\mu,\la)-(\pa+\la+\a_i\mu+\b_i)f_{-1,i}(\pa,\la)\nonumber\\
=&\big((\a_{-1}-1)\mu-\la+\b_{-1}\big)f_{-1,i}(\pa,\la+\mu).
\end{align}
Setting $\mu=0$ in (\ref{eq3.5}), we can obtain that $(\b_{i-1}-\b_{-1}-\b_i)f_{-1,i}(\pa,\la)=0$ for $i\in\Z_{\ge2}$. By (C3) and (\ref{eq3.2}), we have $f_{-1,i}(\pa,\la)\ne0$, which implying $\b_{i-1}=\b_{-1}+\b_i$ for $i\in\Z_{\ge2}$. By a similar discussion, taking (\ref{eq3.3}) into the Jacobi identity $[L_{\mu}[{X_{-1}}_\la X_1]]=[[L_{\mu}X_{-1}]_{\la+\mu} X_1]+[{X_{-1}}_\la [L_{\mu} X_1]]$, we have
\begin{align}\label{eq3.6}
(\pa+&2\mu)g_{-1,1}(\pa+\mu,\la)L+(\pa+\mu)f_{-1,1}(\pa+\mu,\la)H\nonumber\\
=&\big((\a_{-1}-1)\mu-\la+\b_{-1}\big)\big(g_{-1,1}(\pa,\la+\mu)L+f_{-1,1}(\pa,\la+\mu)H\big)\\
&+(\pa+\la+\a_1\mu+\b_1)\big(g_{-1,1}(\pa,\la)L+f_{-1,1}(\pa,\la)H\big)\nonumber.
\end{align}
Letting $\mu=0$ in (\ref{eq3.6}), comparing the coefficients of $L$ and $H$ respectively, we can obtain that $(\b_{-1}+\b_1)g_{-1,1}(\pa,\la)=(\b_{-1}+\b_1)f_{-1,1}(\pa,\la)=0$. By (C3) and (\ref{eq3.3}), we can deduce that $g_{-1,1}(\pa,\la)=0$ and $f_{-1,1}(\pa,\la)=0$ can not hold simultaneously. Thus, we have $\b_{-1}+\b_1=0$. By (\ref{CHV}), we can obtain that $\b_0=0$ in (\ref{eq3.1}). Hence, we have $\b_i=i\b_1$ for $i\in\Z_{\ge-1}$, $\b_1\in\C$.

Similarly, using the Jacobi identity $[{X_{-1}}_{\la}[H_\mu X_i]]=[[{X_{-1}}_{\la}H]_{\la+\mu} X_i]+[H_\mu [{X_{-1}}_{\la} X_i]]$ and comparing the coefficients of $X_{i-1}$ for $i\in\Z_{\ge2}$, we can obtain that
\begin{align}\label{eq3.7}
\g_if_{-1,i}(\pa,\la)=\g_{i-1}f_{-1,i}(\pa+\mu,\la)-\g_{-1}f_{-1,i}(\pa,\la+\mu).
\end{align}
Setting $\mu=0$ in (\ref{eq3.7}), by (C3) and (\ref{eq3.2}), we can deduce that $\g_i=\g_{i-1}-\g_{-1}$ for $i\in\Z_{\ge2}$. By a similar discussion, taking (\ref{eq3.3}) into the Jacobi identity $[H_{\la}[{X_1}_\mu {X_{-1}}]]=[[H_{\la}X_1]_{\la+\mu} {X_{-1}}]+[{X_1}_\mu [H_{\la} {X_{-1}}]]$, we have
\begin{align}\label{eq3.8}
g_{1,-1}&(\pa+\la,\mu)\la H-\g_{-1}\big(g_{1,-1}(\pa,\mu)L+f_{1,-1}(\pa,\mu)H\big)\nonumber\\
&=\g_1\big(g_{1,-1}(\pa,\la+\mu)L+f_{1,-1}(\pa,\la+\mu)H\big).
\end{align}
Comparing the coefficients of $L$ and $H$ in (\ref{eq3.8}) respectively, we can obtain that
\begin{align}
&\g_1g_{1,-1}(\pa,\la+\mu)+\g_{-1}g_{1,-1}(\pa,\mu)=0,\label{eq3.9}\\
&\g_1f_{1,-1}(\pa,\la+\mu)+\g_{-1}f_{1,-1}(\pa,\mu)=g_{1,-1}(\pa+\la,\mu)\la\label{eq3.10}.
\end{align}

From (C3) and (\ref{CHV}), it follows that $\g_{-1}\ne0$ and $\g_0=0$ in (\ref{eq3.1}). If $\g_1=0$, by (\ref{eq3.9}) and (\ref{eq3.10}), we can immediately deduce that $g_{1,-1}(\pa,\mu)=f_{1,-1}(\pa,\mu)=0$, which is a contradiction. Thus, $\g_1\ne0$. It is not difficult to check that the polynomial $g_{1,-1}(\pa,\la)$ is independent of $\la$ in (\ref{eq3.9}), so we abbreviate it to $g_{1,-1}(\pa)$. Then, the formula (\ref{eq3.9}) turns into
$(\g_1+\g_{-1})g_{1,-1}(\pa)=0.$ If $\g_1+\g_{-1}\ne0$, we have $g_{1,-1}(\pa)=0$. Taking this into (\ref{eq3.10}) and by a similar discussion, we can deduce that $f_{1,-1}(\pa,\la)=0$, which is also a contradiction. Then, $\g_1+\g_{-1}=0$. Thus, we have $\g_i=i\g_1$ for $i\in\Z_{\ge-1}$, $\g_1\in\C^*$.

This completes the proof.
\QED
 \vspace{.5cm}

In order to characterize the structure of the Lie conformal algebra ${\A}$, we still need to determine all $\a_i\in\C$ and $g_{-1,1}(\pa,\la),f_{i,j}(\pa,\la)\in \C[\pa,\la]$ for $i,j\in\Z_{\ge-1}$. Let us compute $g_{-1,1}(\pa,\la)$ and $f_{-1,1}(\pa,\la)$ first. Suppose that $f_{1,-1}(\pa,\la)=\sum_{i=0}^{m}a_{i}(\pa)\la^i$ is the solution of (\ref{eq3.10}), where $a_{i}(\pa) \in \C[\pa]$ and $m$ is the highest degree of $\la$ with $a_{m}(\pa)\ne0$.
\vspace{.5cm}

\textbf{Case 1.} $m\ge2$.\\

Taking $f_{1,-1}(\pa,\la)=\sum_{i=0}^{m}a_{i}(\pa)\la^i$ into (\ref{eq3.10}), inspired by the proof of Lemma \ref{lm3.1}, we can obtain that
\begin{align}
g_{1,-1}(\pa+\la)\la =\g_1\sum_{i=0}^{m}a_{i}(\pa)\big((\la+\mu)^i-\mu^i\big).\label{eq3.11}
\end{align}
Comparing the coefficients of $\mu^{m-1}$ in (\ref{eq3.11}), we get $ma_{m}(\pa)\la=0$, which is impossible. Thus, $m\le1$.
 \vspace{.5cm}

\textbf{Case 2.} $m=1$.\\

In this case, we have $f_{1,-1}(\pa,\la)=a_{1}(\pa)\la+a_{0}(\pa)$. It follows from (\ref{eq3.10}) that $g_{1,-1}(\pa+\la)\la =\g_1a_1(\pa)\la$, where $a_1(\pa)\ne0$. Thus, $g_{1,-1}(\pa)$ and $a_1(\pa)$ are nonzero complex numbers. Now, we can suppose that $g_{1,-1}(\pa,\la)=b_{1,-1}^0$ with $b_{1,-1}^0\in \C^*$. Hence, we have $a_1(\pa)=\frac{b_{1,-1}^0}{\g_1}$, i.e., $f_{1,-1}(\pa,\la)=\frac{b_{1,-1}^0}{\g_1}\la+a_{0}(\pa)$.

Letting $\pa=0$ in (\ref{eq3.6}), we can obtain that
\begin{align}\label{eq3.12}
\big(2\mu &g_{-1,1}(\mu,\la)L+\mu f_{-1,1}(\mu,\la)H\big)-(\la+\a_1\mu+\b_1)\big(g_{-1,1}(0,\la)L+f_{-1,1}(0,\la)H\big)\nonumber\\
=&\big((\a_{-1}-1)\mu-\la+\b_{-1}\big)\big(g_{-1,1}(0,\la+\mu)L+f_{-1,1}(0,\la+\mu)H\big).
\end{align}
By the skew-symmetry of $f_{-1,1}(\pa,\la)$ and $g_{-1,1}(\pa)$, we have $f_{-1,1}(\pa,\la)=-f_{1,-1}(\pa,-\pa-\la)=\frac{b_{1,-1}^0}{\g_1}(\pa+\la)-a_{0}(\pa)$ and $g_{-1,1}(\pa,\la)=-g_{1,-1}(\pa,-\pa-\la)=-b_{1,-1}^0$. Taking these two equations into (\ref{eq3.12}) and comparing the coefficients of $L$ and $H$ respectively, we can obtain that
\begin{align}
2\mu b_{1,-1}^0=&b_{1,-1}^0\big((\a_{-1}-1)\mu-\la+\b_{-1}\big)+b_{1,-1}^0(\la+\a_1\mu+\b_1),\label{eq3.13}\\
\mu\Big(\frac{b_{1,-1}^0}{\g_1}(\mu+\la)-a_0(\mu)\Big)=&\big((\a_{-1}-1)\mu-\la+\b_{-1}\big)\Big(\frac{b_{1,-1}^0}{\g_1}(\la+\mu)-a_0(0)\Big)\nonumber\\
&+(\la+\a_1\mu+\b_1)\Big(\frac{b_{1,-1}^0}{\g_1}\la-a_0(0)\Big)\label{eq3.14}.
\end{align}
Since $b_{1,-1}^0\ne0$ and $\b_1+\b_{-1}=0$, by dividing $b_{1,-1}^0$ on both sides of the formula (\ref{eq3.13}) and comparing the coefficients of $\mu$, we can deduce that $\a_1+\a_{-1}=3$. Taking $\a_1+\a_{-1}=3$ and $\b_1+\b_{-1}=0$ into (\ref{eq3.14}), we have
\begin{align*}
\mu\Big(\frac{b_{1,-1}^0}{\g_1}(\mu+\la)-a_0(\mu)\Big)=&\frac{b_{1,-1}^0}{\g_1}\mu\big((\a_{-1}-1)\mu-\la+\b_{-1}\big)+2\mu\Big(\frac{b_{1,-1}^0}{\g_1}\la-a_0(0)\Big).
\end{align*}
Then by direct computation, we get $a_0(\mu)=\frac{b_{1,-1}^0}{\g_1}\big((2-\a_{-1})\mu+\b_{-1}\big)$. Thus, $f_{-1,1}(\pa,\la)=\frac{b_{1,-1}^0}{\g_1}\big((\a_{-1}-1)\pa+\la-\b_{-1}\big)$.

Using the Jacobi identity $[{X_{-1}}_{\la}[{X_{-1}}_\mu X_1]]=[[{X_{-1}}_{\la}{X_{-1}}]_{\la+\mu} X_1]+[{X_{-1}}_\mu [{X_{-1}}_{\la} X_1]]$ and comparing the coefficients of $X_{-1}$, we can obtain that
\begin{align}
&b_{1,-1}^0\big((\a_{-1}-1)\pa+\a_{-1}\la-\b_{-1}\big)+\frac{b_{1,-1}^0}{\g_1}\big((\a_{-1}-1)(\pa+\la)+\mu-\b_{-1}\big)\g_{-1}\nonumber\\
=&b_{1,-1}^0\big((\a_{-1}-1)\pa+\a_{-1}\mu-\b_{-1}\big)+\frac{b_{1,-1}^0}{\g_1}\big((\a_{-1}-1)(\pa+\mu)+\la-\b_{-1}\big)\g_{-1}.\label{eq3.15}
\end{align}
Taking $\g_1+\g_{-1}=0$ into (\ref{eq3.15}), since $b_{1,-1}^0\ne0$, we can deduce that $\la=\mu$, which is impossible. Hence, $m=0$.
 \vspace{.5cm}

\textbf{Case 3.} $m=0$.\\

Since $f_{1,-1}(\pa,\la)=a_{1,-1}^0$ and $\g_1+\g_{-1}=0$, by (\ref{eq3.10}), we obtain $g_{1,-1}(\pa,\la)=0$. By the skew-symmetry of $f_{-1,1}(\pa,\la)$ and $g_{-1,1}(\pa)$, we have $f_{-1,1}(\pa,\la)=-f_{1,-1}(\pa,-\pa-\la)=-a_{1,-1}^0=a_{-1,1}^0$ and $g_{-1,1}(\pa,\la)=-g_{1,-1}(\pa,-\pa-\la)=0$. Then (C3) and (\ref{eq3.3}) lead to $a_{-1,1}^0\in\C^*$.

By the above discussion, we can obtain the following lemma immediately.
\begin{lemm}\label{lm3.2}
$f_{-1,1}(\pa,\la)=a_{-1,1}^0$ and $g_{-1,1}(\pa,\la)=0$, where $a_{-1,1}^0\in\C^*$. In particular, Lie conformal algebra ${\A}$ is not simple.
\end{lemm}
{\it Proof.} It is not difficult to verify that the $\C[\pa]$-submodule of $\A$, which is $\C[\pa]$-spanned by $\{X_i\ |\ i\in\Z_{\ge-1}\}$, is a nontrivial ideal of ${\A}$.
\QED
 \vspace{.5cm}

 Next, we need to determine $f_{-1,i}(\pa,\la), f_{1,i}(\pa,\la)$ for $i\in\Z_{\ge1}$, which are polynomials of $\pa$ and $\la$. Letting $i=2$ in (\ref{eq3.7}), by Lemma \ref{lm3.1}, we get
\begin{align}\label{eq3.16}
2f_{-1,2}(\pa,\la)=f_{-1,2}(\pa,\la+\mu)+f_{-1,2}(\pa+\mu,\la).
\end{align}
Suppose that $f_{-1,2}(\pa,\la)=\sum_{i=0}^{n}a_{i}(\pa)\la^i$ is the solution of (\ref{eq3.16}), where $a_{i}(\pa) \in \C[\pa]$ and $n$ is the highest degree of $\la$ with $a_{n}(\pa)\ne0$.
\vspace{.5cm}

\textbf{Case 1.} $n\ge2$.\\

Taking $f_{-1,2}(\pa,\la)=\sum_{i=0}^{n}a_{i}(\pa)\la^i$ into (\ref{eq3.16}) and comparing the coefficients of $\la^n, \la^{n-1}$ and $\la^{n-2}$ respectively, we can obtain that
\begin{align}
2a_n(\pa)&=a_n(\pa+\mu)+a_n(\pa),\label{eq3.17}\\
2a_{n-1}(\pa)&=a_{n-1}(\pa+\mu)+na_n(\pa)\mu+a_{n-1}(\pa),\label{eq3.18}\\
2a_{n-2}(\pa)&=a_{n-2}(\pa+\mu)+(n-1)a_{n-1}(\pa)\mu+\frac{n(n-1)}{2}a_{n}(\pa)\mu^2+a_{n-2}(\pa).\label{eq3.19}
\end{align}

 One can easily deduce that
 \begin{align}\label{eq3.20}
 a_n(\pa)=k, \quad a_{n-1}(\pa)=p\pa+q, \quad a_{n-2}(\pa)=x\pa^2+y\pa+z
 \end{align}
 from (\ref{eq3.17})--(\ref{eq3.19}), where $k,p,x\in\C^*$, $q,y,z\in\C$. Taking (\ref{eq3.20}) into (\ref{eq3.17})--(\ref{eq3.19}) gives that $p\mu+nk\mu=0$ and $-2x\pa\mu-x\mu^2-y\mu=(p\pa+q)(n-1)\mu+\frac{n(n-1)}{2}k\mu^2$. Now comparing the coefficients of $\pa\mu$ and $\mu^2$ respectively, we get $2x=n(n-1)k$ and $x=-\frac{n(n-1)}{2}k$. That leads to $k=0$, which is a contradiction. Thus, $n\le1$.
  \vspace{.5cm}

\textbf{Case 2.} $n\le1$.\\

Setting $f_{-1,2}(\pa,\la)=a_{1}(\pa)\la+a_{0}(\pa)$, by (\ref{eq3.16}), we can obtain that $f_{-1,2}(\pa,\la)=a(\pa-\la)+b$, where $a,b\in \C$. Letting $i=2$ in (\ref{eq3.5}), we can deduce that
 \begin{align}\label{eq3.21}
(\pa&+\a_{1}\mu+\b_{1})f_{-1,2}(\pa+\mu,\la)-\big((\a_{-1}-1)\mu-\la+\b_{-1}\big)f_{-1,2}(\pa,\la+\mu)\nonumber\\
&=(\pa+\la+\a_2\mu+\b_2)f_{-1,2}(\pa,\la).
\end{align}
 The above formula together with (\ref{eq3.16}) gives
 \begin{align}\label{eq3.22}
 \big(\pa-\la+(\a_2+2\a_{-1}-2)\mu\big)f_{-1,2}(\pa,\la+\mu)=\big(\pa-\la+(2\a_{1}-\a_2)\mu\big)f_{-1,2}(\pa+\mu,\la).
 \end{align}

 By Lemma \ref{lm3.2}, (\ref{eq3.6}) becomes the following formula
 \begin{align*}
(\pa+\mu)a_{-1,1}^0H=\big((\a_{-1}-1)\mu-\la+\b_{-1}\big)a_{-1,1}^0H+(\pa+\la+\a_1\mu+\b_1)a_{-1,1}^0H.
\end{align*}
Comparing the coefficients of $a_{-1,1}^0\mu H$, we obtain $\a_{-1}+\a_{1}=2$. Using $\a_{-1}+\a_{1}=2$ and $f_{-1,2}(\pa,\la)=a(\pa-\la)+b$ in (\ref{eq3.22}), we have
  \begin{align*}
 \big(\pa-\la+(\a_2-2\a_{1}+2)\mu\big)\big(a(\pa-\la-\mu)+b\big)=\big(\pa-\la+(2\a_{1}-\a_2)\mu\big)\big(a(\pa+\mu-\la)+b\big).
 \end{align*}
It follows that $a=0$ and $\a_2-2\a_{1}+1=0$. Thus, $n=0$. By (C3), we get $f_{-1,2}(\pa,\la)=a_{-1,2}^0$, where $a_{-1,2}^0\in\C^*$. Now, we can deduce the following lemma.

 \begin{lemm}\label{lm3.3}
 \emph{(1)} $f_{-1,i}(\pa,\la)=a_{-1,i}^0$ and $f_{1,i}(\pa,\la)=\frac{a_{-1,1}^0}{a_{-1,i+1}^0}\big(\frac{i(i+1)}{2}-1\big)\g_1$, where $a_{-1,i}^0,\g_1\in\C^*$ for $i\in\Z_{\ge1}$.\\
 \emph{(2)} $\a_i=i\a_1-i+1$ for $i\in\Z_{\ge-1}$, where $\a_1\in\C$.
\end{lemm}
{\it Proof.} (1) Using the Jacobi identity $[{X_{-1}}_{\la}[{X_{1}}_\mu X_i]]=[[{X_{-1}}_{\la}{X_{1}}]_{\la+\mu} X_i]+[{X_{1}}_\mu [{X_{-1}}_{\la} X_i]]$ and comparing the coefficients of $X_{i}$, we can obtain that
\begin{align}\label{eq3.23}
f_{1,i}(\pa+\la,\mu)f_{-1,i+1}(\pa,\la)=f_{-1,1}(-\la-\mu,\la)\g_i+f_{-1,i}(\pa+\mu,\la)f_{1,i-1}(\pa,\mu),
\end{align}
for $i\in\Z_{\ge1}$.
Taking $i=1$ in (\ref{eq3.23}), one can easily verify that $a_{-1,2}^0f_{1,1}(\pa+\la,\mu)=a_{-1,1}^0\g_1-a_{-1,1}^0\g_1=0$, which implies $f_{1,1}(\pa,\la)=0$ since $a_{-1,2}^0\ne 0$.

From $f_{-1,2}(\pa,\la)=a_{-1,2}^0$ and $f_{1,1}(\pa,\la)=0$, by (\ref{eq3.23}), we can inductively deduce that
\begin{align}
&f_{-1,i}(\pa,\la)=a_{-1,i}^0,\\
&f_{1,i}(\pa,\la)=\frac{a_{-1,1}^0}{a_{-1,i+1}^0}\big(\frac{i(i+1)}{2}-1\big)\g_1,
\end{align}
where $a_{-1,i}^0,\g_1\in\C^*$ for $i\in\Z_{\ge1}$.

(2) Together with Lemma \ref{lm3.1} and (1), we can obtain that $\a_{i-1}=\a_i+\a_{-1}-1$ by comparing the coefficients of $a_{-1,i}^0\mu X_{i-1}$ for $i\in\Z_{\ge2}$ in (\ref{eq3.5}). Since $\a_{-1}+\a_{1}=2$, $\a_2-2\a_{1}+1=0$ and $\a_0=1$, by (\ref{CHV}), we can inductively deduce that $\a_i=i\a_1-i+1$ for $i\in\Z_{\ge-1}$.
\QED
 \vspace{.5cm}

Finally, we can determine $f_{i,j}(\pa,\la)$ for $i,j\in\Z_{\ge1}$, which are polynomials of $\pa$ and $\la$. By the above discussion, we get the following lemma.
 \begin{lemm}\label{lm3.4}
$f_{i,j}(\pa,\la)=\frac{a_{-1,1}^0a_{-1,2}^0\cdots a_{-1,j}^0}{a_{-1,i+1}^0a_{-1,i+2}^0\cdots a_{-1,i+j}^0}\frac{(i+j+1)!}{(i+1)!(j+1)!}(j-i)\g_1$, where $a_{-1,i}^0,\g_1\in\C^*$ for $i,j\in\Z_{\ge1}$.
\end{lemm}
{\it Proof.} Using the Jacobi identity $[{X_i}_{\la}[{X_{-1}}_\mu X_2]]=[[{X_i}_{\la}{X_{-1}}]_{\la+\mu} X_2]+[{X_{-1}}_\mu [{X_i}_{\la} X_2]]$ and comparing the coefficients of $X_{i+1}$, we can obtain that
\begin{align}\label{eq3.26}
f_{-1,2}(\pa+\la,\mu)f_{i,1}(\pa,\la)=f_{i,-1}(-\la-\mu,\la)f_{i-1,2}(\pa,\la+\mu)+f_{i,2}(\pa+\mu,\la)f_{-1,i+2}(\pa,\mu),
\end{align}
for $i\in\Z_{\ge1}$. By Lemma \ref{lm3.3}, we have $f_{1,2}(\pa,\la)=\frac{2a_{-1,1}^0}{a_{-1,3}^0}\g_1$. Setting $i=2$ in (\ref{eq3.26}), we get $f_{2,2}(\pa,\la)=0$ immediately by Lemma \ref{lm3.3}. From this, using (\ref{eq3.26}) and the skew-symmetry, we can inductively deduce that
\begin{align}\label{eq3.27}
f_{i,2}(\pa,\la)=\frac{a_{-1,1}^0a_{-1,2}^0}{a_{-1,i+1}^0a_{-1,i+2}^0}\frac{(i+3)!}{(i+1)!3!}(2-i)\g_1,
\end{align}
 for $i\in\Z_{\ge1}$.

 Similarly, using the Jacobi identity $[{X_i}_{\la}[{X_{1}}_\mu X_j]]=[[{X_i}_{\la}{X_{1}}]_{\la+\mu} X_j]+[{X_{1}}_\mu [{X_i}_{\la} X_j]]$ and comparing the coefficients of $X_{i+j+1}$, we can obtain that
\begin{align}\label{eq3.28}
f_{i,j+1}(\pa,\la)f_{1,j}(\pa+\la,\mu)=f_{i,1}(-\la-\mu,\la)f_{i+1,j}(\pa,\la+\mu)+f_{i,j}(\pa+\mu,\la)f_{1,i+j}(\pa,\mu),
\end{align}
for $i,j\in\Z_{\ge1}$. Letting $j=2$ in (\ref{eq3.28}), by (\ref{eq3.27}) and Lemma \ref{lm3.3}, we have
\begin{align}\label{eq3.29}
f_{i,3}(\pa,\la)=\frac{a_{-1,1}^0a_{-1,2}^0a_{-1,3}^0}{a_{-1,i+1}^0a_{-1,i+2}^0a_{-1,i+3}^0}\frac{(i+4)!}{(i+1)!4!}(3-i)\g_1,
\end{align}
 for $i\in\Z_{\ge1}$. Then from (\ref{eq3.27})--(\ref{eq3.29}) and Lemma \ref{lm3.3}, we can inductively deduce that
 \begin{align}
 f_{i,j}(\pa,\la)=\frac{a_{-1,1}^0a_{-1,2}^0\cdots a_{-1,j}^0}{a_{-1,i+1}^0a_{-1,i+2}^0\cdots a_{-1,i+j}^0}\frac{(i+j+1)!}{(i+1)!(j+1)!}(j-i)\g_1
 \end{align}
 where $a_{-1,i}^0,\g_1\in\C^*$ for $i,j\in\Z_{\ge1}$.
 \QED
 \vspace{.5cm}

 Now we can prove the main theorem in this section as follows.
 \begin{theo}\label{th3.5}
 Let ${\A}=\oplus^\infty_{i=-1}\A_i$ be a $\Z$-graded conformal algebra satisfying conditions \emph{(C1)--(C3)}. Then $\A\cong{\HV}(\a,\b)$, where $\a,\b\in\C$.
 \end{theo}
 {\it Proof.} Set
\begin{align*}
 L^{'}=L,\quad X_{-1}^{'}={X_{-1}}, \quad X_0^{'}=\frac{1}{\g_1}X_0,\quad X_i^{'}=\frac{(i+1)!}{a_{-1,1}^0a_{-1,2}^0\cdots a_{-1,i}^0\g_1}X_i,
 \end{align*}
 for $i\in\Z_{\ge1}$. It follows from Lemmas \ref{lm3.1}--\ref{lm3.4} and (\ref{eq3.1})--(\ref{eq3.3}) that
\begin{align*}
[{L^{'}}_\la L^{'}]=(\pa+2\la)L^{'}, \quad [{L^{'}}_\la X_i^{'}]=\big(\pa+(i\a_1-i+1)\la+i\b_1\big)X_i^{'}, \quad [{X_i^{'}}_\la X_j^{'}]=(j-i)X_{i+j}^{'},
\end{align*}
for $i,j\in\Z_{\ge-1}$ where $\a_1,\b_1\in\C$. Therefore, $\A\cong{\HV}(\a,\b)$, where $\a,\b\in\C$.

 \QED

\begin{rema}The annihilation algebra of Lie conformal algebra ${\HV}(\a,\b)$ is a Lie algebra $$Lie\big({\HV}(\a,\b)\big)^{+}=span_\C\{L_{n}, H_{i,m}\ |\ L,H_i\in {\HV}(\a,\b), i,n \in \Z_{\ge-1}, m\in\mathbb{Z_{+}}\}$$ with Lie brackets
 \begin{align*}
&[L_m,L_n]=(m-n)L_{m+n}, \quad [H_{i,m}, H_{j,n}]=(j-i)H_{i+j,m+n},\\
&[L_m,H_{i,n}]=\big((m+1)(i\a-i)-n\big)H_{i,m+n}+i\b H_{i,m+n+1},
 \end{align*}
for $i,j\in\Z_{\ge-1}$, where $\a,\b\in\C$.
\end{rema}

\section{Finite nontrivial irreducible modules over graded Lie conformal algebras ${\HV}(\a,\b)$ with $\a\ne1$ and $\b\ne0$}
In this section, we study some properties of graded Lie conformal algebras ${\HV}(\a,\b)$ with $\a,\b\in\C$. Furthermore, we will classify all finite nontrivial irreducible modules over graded Lie conformal algebras ${\HV}(\a,\b)$ with $\a\ne1$ and $\b\ne0$.

\begin{defi}
\begin{em}
Let $\A$ be a Lie conformal algebra and $V$ be a finite free $\A$-module. A finite chain of $\A$-submodules of $V$
\begin{align}\label{eq4.1}
0=V_0\subset V_1 \subset V_2\subset\cdots\subset V_n\subseteq V,
\end{align}
is called \emph{free compositions series of nontrivial index $n$} such that $V/V_n$ is a trivial free $\A$-module and $V_i/V_{i-1}$ for $i=1,2,3,\cdots,n$ is a free extension of a nontrivial irreducible $\A$-module and a torsion $\C[\pa]$-module(i.e., a trivial $\A$-module).

In particular, if $V_n=V$, then $V$ is called \emph{completely nontrivial}. We use $l_V$ to denote the length of free compositions series of $V$, which is minimal number of all nontrivial index of free compositions series of $V$.
\end{em}
\end{defi}

Now, let us introduce the definition of free artianian Lie conformal algebra.
\begin{defi}
\begin{em}
A Lie conformal algebra $\A$ is \emph{free artianian} if any finite nontrivial free $\A$-module always has a nontrivial irreducible $\A$-submodule.
\end{em}
\end{defi}

Inspired by \cite{CK}, one can obtain that any finite semisimple Lie conformal algebra is free artianian. In \cite{XH}, some properties of free artianian Lie conformal algebra are investigated.
\begin{prop}\label{pro4.3}\emph{(\cite{XH}, Proposition 3.1)}
Let $\A$ be a free artianian Lie conformal algebra. Then any finite free $\A$-module has a free compositions series.
\end{prop}

Let $\mathcal{M}$ be a free $\C[\pa]$-module with a $\C[\pa]$-basis $\{J_{(a,b)}\ |\ (a,b)\in\C\times\C\}$. Then $\mathcal{M}$ is a $Vir$-module by defining
\begin{align}
L_\la J_{(a,b)}=(\pa+a\la+b)J_{(a,b)},\quad \forall\ (a,b)\in\C\times\C.
\end{align}
\begin{prop}\emph{(\cite{XH}, Proposition 3.2)}
Let V be a finite nontrivial free $Vir\ltimes \mathcal{M}$-module. Then there are only finitely many $(a,b)\in\C\times\C$ such that $J_{(a,b)}$ acts nontrivially on V.
\end{prop}

\begin{lemm}\label{lem4.4}
Let $V$ be a finite nontrivial free ${\HV}(\a,\b)$-module with $\a\ne1$ and $\b\ne0$. Then there are only finitely many $H_i\in{\HV}(\a,\b)$ such that $H_i$ acts nontrivially on $V$, where $H_i$ as shown in Definition \ref{def2.7}.
\end{lemm}
 {\it Proof.} If $L$ acts trivially on $V$, then for any $i\in\Z_{\ge-1}$ and $v\in V$, we have $[L_\la H_i]_{\la+\mu} v=\big((i\a-i)\la-\mu+i\b\big){H_i}_{\la+\mu}v=0$, which implies ${H_i}_{\la}v=0$ for $i\in\Z_{\ge-1}$. This is a contradiction with that $V$ is a nontrivial ${\HV}(\a,\b)$-module. Thus, $L$ acts nontrivially on $V$.

 Since $Vir$ is artianian, by Proposition \ref{pro4.3}, $V$ has a free compositions series of nontrivial index $n$ such that $V/V_n$ is a trivial free $Vir$-module and $V_i/V_{i-1}$ is nontrivial irreducible $Vir$-module for $i=1,2,3,\cdots,n$. Suppose that $V_n$ has a $\C[\pa]$-basis $\{v_1,v_2,\cdots,v_n\}$ such that each $V_i/V_{i-1}=\C[\pa]\bar{v_i}$ is a rank one free $Vir$-module $V_{a_i,b_i}$ where $a_i\ne0$ for $i=1,2,3,\cdots,n$.
   \vspace{.5cm}

  \textbf{Case 1.} ${H_i}_\la V_n\nsubseteq V_n[\la]$.\\

Suppose that $\{\bar{u_1},\bar{u_2},\cdots,\bar{u_t}\}$ be minimal generators of $V/V_n$ as a $\C[\pa]$-module. For each $j$, let $u_j$ be the preimage of $\bar{u_j}$ in $V$. Assume $k$ is the smallest integer such that ${H_i}_\la v_k=\sum_jg_j(\pa,\la)u_j+t$ for some $t\in V_n[\la]$ and $0\ne\sum_jg_j(\pa,\la)u_j\notin V_n[\la]$. Applying $L_\la$ to ${H_i}_\mu v_k$, we have
\begin{align}
 \big(-\mu+(i\a-i)\la+i\b\big)\sum_jg_j(\pa,\la+\mu)u_j=-(\pa+\mu+a_k\la+b_k)\sum_jg_j(\pa,\mu)u_j \ \  \text{mod}\ V_n[\la].
\end{align}
Comparing the degree of $\pa$, we obtain a contradiction. Hence, we must have ${H_i}_\la V_n\subseteq V_n[\la]$.\\
   \vspace{.5cm}

 \textbf{Case 2.} ${H_i}_\la V_n\subseteq V_n[\la]$.\\

Since $V/V_n$ is a trivial $Vir$-module, if $H_i$ acts trivially on $V_n$, we have
 \begin{align*}
 {\big(-\mu+(i\a-i)\la+i\b\big)H_i}_{\la+\mu}v=[L_\la H_i]_{\la+\mu}v=L_\la({H_i}_\mu v)-{H_i}_\mu(L_\la v)=L_\la({H_i}_\mu v),
 \end{align*}
for any $v\in V$. Comparing the degree of $\mu$ leads to ${H_i}_\mu v=0$. Thus, we suppose that $H_i$ acts nontrivially on $V_n$ for $i\in\Z^*_{\ge-1}$. For each $j$, there exists some $k$ such that $${H_i}_\mu v_j=f_{jk}(\pa,\mu)v_k\quad \text{mod}\ V_{k-1}[\mu],$$ where $f_{jk}(\pa,\mu)\ne0$. Applying $L_\la$ to ${H_i}_\mu v_j$ and comparing the coefficients of $v_k$, we have
\begin{align}\label{eq4.3}
(\pa&+a_k\la+b_k)f_{jk}(\pa+\la,\mu)-(\pa+\mu+a_j\la+b_j)f_{jk}(\pa,\mu) \nonumber\\
 &=\big(-\mu+(i\a-i)\la+i\b\big)f_{jk}(\pa,\la+\mu).
\end{align}
Letting $\la=0$ in (\ref{eq4.3}), one can easily verify that $i\b=b_k-b_j$. Since $b_k,b_j$ have only finitely many choices, the choices of $i\b$ are finite. Suppose that the total degree of $f_{jk}(\pa,\la)$ is $m$. Comparing the terms of degree $m+1$ in (\ref{eq4.3}), we obtain
\begin{align}\label{eq4.4}
(\pa+a_k\la)f_{jk}^m(\pa+\la,\mu)-(\pa+\mu+a_j\la)f_{jk}^m(\pa,\mu)= \big(-\mu+(i\a-i)\la\big)f_{jk}^m(\pa,\la+\mu),
\end{align}
where $f_{jk}^m(\pa,\la)\ne0$ is the $m$-th homogeneous part of $f_{jk}(\pa,\la)$. Setting $\mu=0$ in (\ref{eq4.4}), it follows that
\begin{align}\label{eq4.5}
(\pa+a_k\la)f_{jk}^m(\pa+\la,0)-(\pa+a_j\la)f_{jk}^m(\pa,0)= (i\a-i)\la f_{jk}^m(\pa,\la).
\end{align}
Since $i\in\Z^*_{\ge-1}$ and $\a\ne1$, we have $f_{jk}^m(\pa,0)\ne0$. Comparing the coefficients of $\la\pa^m$ in (\ref{eq4.5}), we get
\begin{align}\label{eq4.6}
m=(i\a-i)+a_j-a_k.
\end{align}
Taking $m=(i\a-i)+a_j-a_k$ in (\ref{eq4.4}), for each $s,t\in\{0,1,2,\cdots,m+1\}$, the coefficients of $\pa^s\la^t\mu^{m+1-s-t}$ in (\ref{eq4.4}) provide a series of polynomial equations of $i\a-i+1$ with coefficients associated with $a_j$ and $a_k$. Note that the choices of $a_j$ and $a_k$ are finite. Since these polynomial equations are only determined by $a_j, a_k$ and (\ref{eq4.4}), the choices of different $i\a-i+1$ are also finite.

Since $\a\ne1$ and $\b\ne0$, by the above discussion, we can deduce that there are only finitely many $H_i\in{\HV}(\a,\b)$ such that $H_i$ acts nontrivially on $V$.
 \QED

\begin{coro}\label{cor4.5}
For $\a\ne1$ and $\b\ne0$, the graded Lie conformal algebra ${\HV}(\a,\b)$ does not have nontrivial finite faithful modules.
\end{coro}
 {\it Proof.} The conclusion is followed by Lemma \ref{lem4.4} and Definition \ref{def2.11}.
 \QED

\begin{lemm}\label{lem4.6}
The set of locally nilpotent elements of the graded Lie conformal algebra ${\HV}(\a,\b)$ with $\alpha,\beta \in \mathbb{C}$ is equal to $\C[\pa]H_{-1}$.
 \end{lemm}
 {\it Proof.} Denote by $\mathcal{N}$ the set of locally nilpotent elements of ${\HV}(\a,\b)$. By (\ref{eq2.5}) and conformal sesquilinearity, we can obtain that $[{f_{-1}(\pa)H_{-1}}_{\la}f_i(\pa)H_{i}]=(i+1)f_{-1}(-\la)f_i(\pa+\la)H_{i-1}$ and $[{f_{-1}(\pa)H_{-1}}_{\la}{g_0(\pa)L}]=f_{-1}(-\la)g_0(\pa+\la)\big((1-\a)\pa+(2-\a)\la+\b\big)H_{-1}$, where $f_{-1}(\pa), g_0(\pa), f_i(\pa)\in\C[\pa]$ for $i\in \Z_{+}$. Then by  (\ref{eq2.5}),  we have $\C[\pa]H_{-1}\subset\mathcal{N}$.

 On the other hand, by (\ref{eq2.5}), it is not difficult to verify $\C[\pa]L\notin \mathcal{N}$. Then we can suppose that $x=\sum_{i=-1}^nf_i(\pa)H_i\in\mathcal{N}$, where $f_i(\pa)\in\C[\pa]$ for $i=-1,0,1,2,...,n$ and $n$ is the maximal integer such that $f_n(\pa)\ne0$. If $n\ge0$, by (\ref{eq2.5}), we get $[{H_n}_\la H_j]=(j-n)H_{n+j}$ for $j\in\Z_+$. For $j\ne n$ and $j\in\Z_+$, applying $(ad\  x)_\la^m$ to $H_j$, we can deduce that the coefficients of $H_{mn+j}$ in the expression of $(ad\ x)_\la^m(H_j)$ is nonzero for $m\gg0$. This is a contradiction with $x\in\mathcal{N}$. Thus, we must have $n=-1$ and the lemma follows.

 \QED

Next, let us investigate all free nontrivial ${\HV}(\a,\b)$-modules of rank one over $\C[\pa]$.
  \begin{proposition}\label{pro4.9}
 All free nontrivial ${\HV}(\a,\b)$-modules of rank one over $\mathbb{C}[\partial]$ with $\alpha,\beta \in \mathbb{C}$ are as follows:
  \begin{align}\label{MHVab}
V_{a,b}=\mathbb{C}[\partial]v,\qquad L_\lambda v=(\partial+a\lambda+b)v,\quad {H_i}_\lambda v=0,
 \end{align}
 where $a,b\in \mathbb{C}$ and $i\in\Z_{\ge-1}$.
 \end{proposition}
 {\it Proof.} Suppose that $V=\mathbb{C}[\partial]v$ is a nontrivial ${\HV}(\a,\b)$-module with $\alpha,\beta \in \mathbb{C}$. Similar to the proof of Lemma \ref{lem4.4}, it is not difficult to verify that $L$ acts nontrivially on $V$. Since ${\HV}(\a,\b)_0=\C[\pa]L\oplus\C[\pa]H_0\cong\HV$, by Proposition \ref{pr2}, we can deduce that $L_\lambda v=(\partial+a\lambda+b)v$ and ${H_0}_\lambda v=c v$ for $a,b,c\in \mathbb{C}$. Assume that ${H_i}_\lambda v=f_i(\pa,\la)v$ for $i\in\Z^*_{\ge-1}$. Using $[{H_0}_{\la}{H_{i}}]_{\la+\mu} v={H_0}_\la{H_{i}}_\mu v-{H_i}_\mu{H_{0}}_\la v$, we can obtain that
 \begin{align}\label{eq4.9}
 if_i(\pa,\la+\mu)v=c\big(f_i(\pa+\la,\mu)-f_i(\pa,\mu)\big)v,
 \end{align}
for $i\in\Z^*_{\ge-1}$.

If $c=0$, by (\ref{eq4.9}), we can immediately deduce that $f_i(\pa,\la)=0$ for $i\in\Z^*_{\ge-1}$.

If $c\ne0$, from $[L_{\la}{H_{i}}]_{\la+\mu} v=L_\la{H_{i}}_\mu v-{H_i}_\mu{L}_\la v$, it follows that
 \begin{align}\label{eq4.10}
 \big((i\a-i)\la-\mu+i\b\big)f_i(\pa,\la+\mu)v=(\pa+a\la+b)f_i(\pa+\la,\mu)v-(\pa+\mu+a\la+b)f_i(\pa,\mu)v,
 \end{align}
 for $i\in\Z^*_{\ge-1}$. Taking (\ref{eq4.9}) into (\ref{eq4.10}), we get
  \begin{align}\label{eq4.11}
 \big((i\a-i)\la-\mu+i\b\big)f_i(\pa,\la+\mu)v=\frac{i}{c}(\pa+a\la+b)f_i(\pa,\la+\mu)v-\mu f_i(\pa,\mu)v,
 \end{align}
 for $i\in\Z^*_{\ge-1}$. Comparing the degree of $\pa$ in (\ref{eq4.11}), we can deduce that $f_i(\pa,\la)=0$ for $i\in\Z^*_{\ge-1}$. However, by $[{H_{-1}}_{\la}{H_{1}}]_{\la+\mu} v={H_{-1}}_\la{H_{1}}_\mu v-{H_1}_\mu{H_{-1}}_\la v$, we have $2{H_0}_{\la+\mu}v=0$, i.e., $2cv=0$, which implies $c=0$.
 This is a contradiction.

 \QED

 \begin{proposition}\label{pro4.10}
 For $\a,\b,a,b\in\C$, $V_{a,b}$ is an irreducible ${\HV}(\a,\b)$-module if and only if $a\ne0$.
\end{proposition}
 {\it Proof.}
In this case, ${H_i}_\la v=0$ for all $i\in\Z_{\ge-1}$. The irreducibility of $V_{a,b}$ as an ${\HV}(\a,\b)$-module is equivalent to that as a $Vir$-module. Therefore, by Proposition \ref{pr1}, we obtain the conclusion.
 \QED

By the above discussion, we can investigate our main result in this section.
 \begin{theo}\label{th4.11}
For $\a\ne1$ and $\b\ne0$, any finite nontrivial irreducible ${\HV}(\a,\b)$-module $V$ is free of rank one over $\C[\pa]$. In particular, $V\cong V_{a,b}$ with $a\ne0$ as shown in \emph{(\ref{MHVab})}.
 \end{theo}
 {\it Proof.} Let $V$ be a finite nontrivial irreducible ${\HV}(\a,\b)$-module. By Lemma \ref{lem4.4}, we can suppose that $i\in\Z_{\ge-1}$ is the biggest integer such that $H_i$ acts nontrivially on $V$. For any $v\in V$, using $[{H_{-1}}_\la H_{i+1}]_{\la+\mu}v={H_{-1}}_\la {H_{i+1}}_{\mu}v-{H_{i+1}}_\mu {H_{-1}}_{\la}v$, we have $(i+2){H_{i}}_{\la+\mu}v=0$, which is a contradiction. Thus, $H_i$ acts trivially on $V$ for all $i\in\Z_{\ge-1}$.

 Since $\oplus_{i\in\Z_{\ge-1}}\C[\pa]H_i$ is a nontrivial ideal of ${\HV}(\a,\b)$, then $V$ can be regarded as a finite nontrivial irreducible module of ${\HV}(\a,\b)/\oplus_{i\in\Z_{\ge-1}}\C[\pa]H_i\cong Vir$. By Proposition \ref{pr1}, we can deduce that $V$ is free of rank one. Together with Proposition \ref{pro4.9}, the theorem holds.
 \QED
\section{Conformal derivations of graded Lie conformal algebras ${\HV}(\a,\b)$}
In this section, we study the conformal derivations of graded Lie conformal algebras ${\HV}(\a,\b)$ with $\a,\b\in\C$.

\begin{defi}
	\em{
		A \emph {conformal linear map} between $\mathbb{C}[\partial]$-modules $\mathcal{A}$ and $\mathcal{B}$ is a $\mathbb{C}$-linear map $\phi_{\lambda}$: $\mathcal{A} \to \mathbb{C}[\lambda] \otimes \mathcal{B}$, satisfying the following axiom:
		\begin{align*}
		\phi_{\lambda}(\partial a)= (\partial+ \lambda)\phi_{\lambda}(a), \quad \forall \ a\in\mathcal{A}.
		\end{align*}}
			\end{defi}
Obviously, $\phi_{\lambda}$ does not depend on the choice of the indeterminate variable $\lambda$. The vector space of all conformal linear maps from $\mathcal{A}$ to $\mathcal{B}$ is denoted by $Chom(\mathcal{A}, \mathcal{B})$, which can be made into a $\mathbb{C}[\partial]$-module via $(\partial \phi)_{\lambda}(a)= -\lambda\phi_{\lambda}(a),$ for all $a\in\mathcal{A}$. For convenience, we write $Cend{\mathcal{A}}$ for $Chom(\mathcal{A}, \mathcal{A})$.

\begin{defi}
	\em{Let $\mathcal{A}$ be a Lie conformal algebra. A conformal linear map $\phi_{\lambda}$: $\mathcal{A} \to \mathbb{C}[\lambda]\otimes\mathcal{A}$ is a \emph {conformal derivation} of $\mathcal{A}$ if
		\begin{align*}
		\phi_{\lambda}([a_{\mu}b])= [\phi_{\lambda}(a)_{\lambda+ \mu}b]+ [a_{\mu}\phi_{\lambda}(b)], \quad \forall \ a, b \in \mathcal{A}.
		\end{align*}}
		\end{defi}

By the Jacobi identity, it is not difficult to check that for every $a \in \mathcal{A}$, the map $(ad\ a) _{\lambda}$ which is defined by $(ad\ a) _{\lambda} b= [a_{\lambda} b]$ for any $b \in \mathcal{A}$, is a conformal derivation of $\mathcal{A}$. All conformal derivations of this kind are called \emph{inner conformal derivations}. Denote by $CDer(\mathcal{A})$ and $CInn(\mathcal{A})$ the vector spaces of all conformal derivations and inner conformal derivations of $\mathcal{A}$ respectively.

There is an important example of a non-inner conformal derivation defined as follows,
\begin{example}\label{ex2.6}
	\em{Let $Cur\mathcal{G}$ be the current Lie conformal algebra associated to the finite dimensional Lie algebra $\mathcal{G}$. Define a conformal linear map $d^L_{\lambda}$: $Cur\mathcal{G} \to Cur\mathcal{G}$ by $d^L_{\lambda}a= (\partial+ \lambda)a$ for every $a \in \mathcal{G}$. One can easily verify that $d^L_{\lambda}$ is a conformal derivation.}

\end{example}

\begin{defi}
	\em{Let $\mathcal{A}$ be a finite Lie conformal algebra. The $\lambda$-bracket on $Cend(\mathcal{A})$ is given by
	\begin{align*}
	[\phi_{\lambda}\psi]_{\mu}a= \phi_{\lambda}(\psi_{\mu- \lambda}a)- \psi_{\mu- \lambda}(\phi_{\lambda}a), \quad \forall \ a \in \mathcal{A},
	\end{align*}
	defines a Lie conformal algebra structure on $Cend(\mathcal{A})$. This is called the \emph{general Lie conformal algebra} on $\mathcal{A}$ which is denoted by $gc(\mathcal{A})$.}

\end{defi}

Inspired by the above definition, it is easy to check that $CDer(\mathcal{A})$ and $CInn(\mathcal{A})$ are Lie conformal subalgebras of $gc(\mathcal{A})$, where $\mathcal{A}$ is a finite Lie conformal algebra. Moreover, we have the following tower $CInn(\mathcal{A})\subseteq CDer(\mathcal{A})\subseteq gc(\mathcal{A})$.

It was also shown in \cite{DK} that conformal derivations on $Vir$ and all simple current Lie conformal algebras $Cur\mathcal{G}$ are as follows.\begin{proposition}
	\emph{(1)} Every conformal derivation of  $Vir$ is inner. \\
	\emph{(2)} For a finite dimensional simple Lie algebra $\mathcal{G}$, every conformal derivation of $Cur\mathcal{G}$ is of the form $p(\partial)d^L_\la+d_\la,$ where $d_\la$ is an inner conformal derivation and $d^L_\la$ is as in Example \ref{ex2.6}.
\end{proposition}

Set $D_\la\in CDer\big({\HV}(\a,\b)\big)$. Define $D^i_\la(X_j)=\pi_{i+j}\big(D_\la(X_j)\big)$, where $\pi_i$ is the natural projection from $\C[\la]\otimes {\HV}(\a,\b)$ onto $\C[\la]\otimes {\HV}(\a,\b)_i$ and $X_i$ is any one of $\C[\partial]$-generators of ${\HV}(\a,\b)_i$ for $i\in \Z_{\ge -1}$. It is not difficult to verify that $D^i_\la$ is a conformal derivation. Moreover, $D_\la=\sum_{i\in \Z_{\ge -1}}D^i_\la$ in the sense that for any $x\in{\HV}(\a,\b)$, there are only finitely many $D^i_\la(x)\ne0$.

\begin{lemm}\label{lem5.6}
For $i\in \Z_{\ge -1}$, $D^i_\la$ is an inner conformal derivation of the form $D^i_\la=\big(ad\ f(\pa)H_i\big)_\la$ for some $f(\pa)\in\C[\pa]$.
\end{lemm}
 {\it Proof.}
 Suppose that
 \begin{align}
 &D^i_\la(H_j)=\pi_{i+j}\big(D_\la(H_j)\big)=f^i_j(\pa,\la)H_{i+j},\label{eq5.1}\\
 &D^i_\la(L)=\pi_{i}\big(D_\la(L)\big)=g^i_0(\pa,\la)H_{i},\label{eq5.2}\\
 &D^1_\la(H_{-1})=\pi_{0}\big(D_\la(H_{-1})\big)=g^1_{-1}(\pa,\la)L+f^1_{-1}(\pa,\la)H_{0},\label{eq5.3}\\
 &D^{-1}_\la(H_1)=\pi_{0}\big(D_\la(H_{1})\big)=g^{-1}_{1}(\pa,\la)L+f^{-1}_{1}(\pa,\la)H_{0},\label{eq5.4}\\
 &D^0_\la(H_k)=\pi_{k}\big(D_\la(H_{k})\big)=f^{0}_{k}(\pa,\la)H_{k},\label{eq5.5}\\
 &D^0_\la(H_0)=\pi_{0}\big(D_\la(H_{0})\big)=g^{0}_{0}(\pa,\la)L+f^{0}_{0}(\pa,\la)H_{0},\label{eq5.6}\\
 &D^0_\la(L)=\pi_{0}\big(D_\la(L)\big)=g^{0}_{01}(\pa,\la)L+f^{0}_{01}(\pa,\la)H_{0},\label{eq5.7}
 \end{align}
for $i,j,k\in\Z_{\ge-1}$ with $i,k\ne0$ and $i+j\ne0$, where $f^n_m(\pa,\la), g^i_0(\pa,\la), g^1_{-1}(\pa,\la), g^{-1}_{1}(\pa,\la), g^{0}_{01}(\pa,\la),\\ f^{0}_{01}(\pa,\la)\in\C[\pa,\la]$ for $i,m,n\in\Z_{\ge-1}$.
  \vspace{.5cm}

\textbf{Case 1.} $i,j\in\Z_{\ge-1}$, $i\ne0$.\\
  \vspace{.5cm}

\textbf{Subcase 1.1.} $i+j\ne0$.\\

Applying $D^i_\la$ to $[{H_0}_\mu H_j]=jH_j$, we can get
\begin{align}\label{eq5.8}
jf^i_j(\pa,\la)=(j-i)f^i_0(-\la-\mu,\la)+(i+j)f^i_j(\pa+\mu,\la).
\end{align}
Setting $\mu=0$ in (\ref{eq5.8}), we can obtain
\begin{align}\label{eq5.9}
jf^i_j(\pa,\la)=(j-i)f^i_0(-\la,\la)+(i+j)f^i_j(\pa,\la).
\end{align}

 If $i=j$, we have $f^i_i(\pa,\la)=0$ by (\ref{eq5.9}).

 If $i\ne j$, we can deduce that $f^i_j(\pa,\la)=\frac{i-j}{i}f^i_0(-\la,\la)$. Similarly, applying $D^i_\la$ to $[L_\mu H_j]=\big(\pa+(j\a-j+1)\mu+j\b\big)H_j$, we can get
\begin{align}\label{eq5.10}
(j-&i)g^i_0(-\la-\mu,\la)+\Big(\pa+\big((i+j)\a-i-j+1\big)\mu+(i+j)\b\Big)f^i_j(\pa+\mu,\la)\nonumber\\
=&\big(\pa+\la+(j\a-j+1)\mu+j\b\big)f^i_j(\pa,\la).
\end{align}
 Plugging $f^i_j(\pa,\la)=\frac{i-j}{i}f^i_0(-\la,\la)$ into (\ref{eq5.10}), we have $g^i_0(-\la-\mu,\la)=\frac{1}{i}f^i_0(-\la,\la)\big((i\a-i)\mu-\la+i\b\big)$. Set $-\la-\mu=\pa$, we can obtain that $g^i_0(\pa,\la)=-\frac{1}{i}f^i_0(-\la,\la)\big(i(\a-1)\pa+(i\a-i+1)\la-i\b\big)$.
   \vspace{.5cm}

\textbf{Subcase 1.2.} $i+j=0$.\\

Applying $D^1_\la$ to $[{H_0}_\mu H_{-1}]=-H_{-1}$, we can get
\begin{align}\label{eq5.11}
g^1_{-1}(\pa,\la)L+f^1_{-1}(\pa,\la)H_0=2f^1_0(-\la-\mu,\la)H_0-g^1_{-1}(\pa+\mu,\la)\mu H_0.
\end{align}
Comparing the coefficients of $L, H_0$, we have $g^1_{-1}(\pa,\la)=0$ and $f^1_{-1}(\pa,\la)=2f^1_0(-\la-\mu,\la)$, which implies that $f^1_{-1}(\pa,\la)$ and $f^1_0(\pa,\la)$ are independent of $\pa$. We abbreviate $f^1_0(\pa,\la)$ to $f^1_0(\la)$. Thus, $f^1_{-1}(\pa,\la)=2f^1_0(\la)$.

Similarly, applying $D^{-1}_\la$ to $[{H_0}_\mu H_{1}]=H_{1}$, we can deduce that  $g^{-1}_{1}(\pa,\la)=0$ and $f^{-1}_{1}(\pa,\la)=2f^{-1}_0(\la)$.

By the above discussion, for any $i\in\Z_{\ge-1}^*$, take $f(\la)=-\frac{1}{i}f^i_0(-\la,\la)$, then we have $D^i_\la (H_j)=\big(ad\ f(-\pa)H_i\big)_\la H_j$ for $j\in\Z_{\ge-1}$.
\vspace{.5cm}

\textbf{Case 2.} $i=0$, $j\in\Z_{\ge-1}$.\\

Applying $D^0_\la$ to $[L_\mu H_{n}]$, $[{L}_\mu H_{0}]$ and $[{H_{0}}_\mu H_{n}]$, respectively, by (\ref{eq5.5})--(\ref{eq5.7}), we can obtain that
\begin{align}
\big(\pa+&\la+(n\a-n+1)\mu+n\b\big)f^0_n(\pa,\la)\nonumber\\
=&nf^0_{01}(-\la-\mu,\la)+\big(\pa+(n\a-n+1)\mu+n\b\big)f^0_n(\pa+\mu,\la)\nonumber\\
&+\big(\pa+(n\a-n+1)(\la+\mu)+n\b\big)g^0_{01}(-\la-\mu,\la)\label{eq5.16},\\
(\pa+&\la+\mu)\big(g^{0}_{0}(\pa,\la)L+f^{0}_{0}(\pa,\la)H_0\big)-(\pa+\la+\mu)g^0_{01}(-\la-\mu,\la)H_0\nonumber\\
=&(\pa+\mu)f^{0}_{0}(\pa+\mu,\la)H_0+(\pa+2\mu)g^{0}_{0}(\pa+\mu,\la)L,\label{eq5.12}\\
nf^0_n&(\pa,\la)-nf^0_n(\pa+\mu,\la)-nf^{0}_{0}(-\la-\mu,\la)\nonumber\\
=&\big(\pa+(n\a-n+1)(\la+\mu)+n\b\big)g^{0}_{0}(-\la-\mu,\la),\label{eq5.14}
\end{align}
for $n\in\Z^*_{\ge-1}$. It is easy to check that $g^{0}_{0}(\pa,\la)=0$ by comparing the degree of $\mu$ in coefficients of $L$ in (\ref{eq5.12}). Then, applying $D^0_\la$ to $[{H_{-1}}_\mu H_{m}]$, we have
\begin{align}
f^0_{m-1}(\pa,\la)=&f^0_{-1}(-\la-\mu,\la)+f^0_{m}(\pa+\mu,\la),\label{eq5.15}
\end{align}
for $m\in\Z_{+}$.

Taking $g^{0}_{0}(\pa,\la)=0$ and (\ref{eq5.15}) with case $m=1$ into (\ref{eq5.14}), we have $f^0_n(\pa,\la)-f^0_n(\pa+\mu,\la)=f^0_{-1}(-\la-\mu,\la)+f^0_1(-\la,\la)$. Setting $-\la-\mu=\pa$ in this formula, we can get
\begin{align}\label{eq5.17}
f^0_n(-\la-\mu,\la)-f^0_n(-\la,\la)=f^0_{-1}(-\la-\mu,\la)+f^0_1(\pa+\mu,\la).
\end{align}
Letting $n=-1$ in (\ref{eq5.17}), we obtain $f^0_1(\pa,\la)=-f^0_{-1}(-\la,\la)$. Then, by (\ref{eq5.15}) with case $m=1$, we can deduce that $f^0_{-1}(\pa,\la)=f^0_{-1}(\la)$ and $f^{0}_{0}(\pa,\la)=f^0_{-1}(\la)+f^0_1(\pa+\mu,\la)=0$. Hence, setting $g^{0}_{0}(\pa,\la)=f^{0}_{0}(\pa,\la)=0$ in (\ref{eq5.12}) and (\ref{eq5.14}) implies $g^0_{01}(\pa,\la)=0$ and $f^0_n(\pa,\la)=f^0_n(\la)$. Together with (\ref{eq5.15}), we can get $f^0_n(\pa,\la)=nf^0_1(\la)$. Setting $g^0_{01}(\pa,\la)=0$ and $f^0_n(\pa,\la)=nf^0_1(\la)$ in (\ref{eq5.16}), we have $f^0_{01}(-\la-\mu,\la)=\la f^0_1(\la)$, i.e., $f^0_{01}(\pa,\la)=\la f^0_1(\la)$.

By the above discussion, take $g(\la)=f^0_1(\la)$, then $D^0_\la (H_j)=\big(ad\ g(-\pa)H_0\big)_\la H_j$ for $j\in\Z_{\ge-1}$ and $D^0_\la (L)=\big(ad\ g(-\pa)H_0\big)_\la L$.

This completes the proof.
 \QED

\begin{theo}\label{th5.7}
Any conformal derivation of graded Lie conformal algebra ${\HV}(\a,\b)$ with $\a,\b\in\C$ is inner, i.e., $CDer\big({\HV}(\a,\b)\big)=CInn\big({\HV}(\a,\b)\big).$
\end{theo}
 {\it Proof.} Suppose that $D_\la=\sum_{i\in \Z_{\ge -1}}D^i_\la$ is a conformal derivation of ${\HV}(\a,\b)$. By Lemma \ref{lem5.6}, we have $D_\la=\sum_{i\in \Z_{\ge -1}}\big(ad\ f_i(\pa)H_i\big)_\la$, where $f_i(\pa)\in\C[\pa]$. If there are infinite many $i$ such that $f_i(\pa)\ne0$, then $D_\la(H_0)=\sum_{i\in \Z_{\ge -1}}-if_i(-\la)H_i$ is an infinite sum, which contradicts with the definition of conformal derivation. Thus, $D_\la=\sum_{i\in \Z_{\ge -1}}\big(ad\ f_i(\pa)H_i\big)_\la$ is a finite sum. Set $x=\sum_{i\in \Z_{\ge -1}}f_i(\pa)H_i\in{\HV}(\a,\b)$. Then $D_\la=(ad\ x)_\la$, which implies $D_\la$ is an inner derivation of ${\HV}(\a,\b)$. Hence, the conclusion holds.

  \QED

\vspace{4mm} \noindent\bf{\footnotesize Acknowledgements.}\ \rm
{\footnotesize This work was supported by the National Natural Science Foundation of China (Nos. 11871421, 11971350, 12171129) and the Zhejiang Provincial Natural Science Foundation of China (No. LY20A010022) and the Fundamental Research Funds for the Central Universities (No. 22120210554).}\\
\vskip18pt \small\footnotesize
\parskip0pt\lineskip1pt

\end{document}